\documentclass{amsart}
\usepackage{amsmath}
\usepackage{amssymb}
\ifx\pdftexversion\undefined
  \usepackage[dvips]{graphics}
\else
  \usepackage[pdftex]{graphics}
\fi
\newtheorem{Proposition}{Proposition}
  \newtheorem{Remark}{Remark}
  \newtheorem{Corollary}[Proposition]{Corollary}
  
   \newtheorem{Theorem}[Proposition]{Theorem}

\def\z{\noindent}
\def\pint{{- \kern-1.1em\int}}
\def\pist#1#2{\noindent\hangindent 2em\hangafter1\hbox to 2em{#1\hfil~~}#2}

\def\Box{{\hfill\hbox{\enspace${\sqre}$}} \smallskip}
\def\sqr#1#2{{\vcenter{\vbox{\hrule height .#2pt
                             \hbox{\vrule width .#2pt height#1pt \kern#1pt
                                   \vrule width .#2pt}
                             \hrule height .#2pt}}}}
\def\sqre{\mathchoice\sqr54\sqr54\sqr{4.1}3\sqr{3.5}3}

\def\CC{\mathbb{C}}

\def\RR{\mathbb{R}}
\def\ZZ{\mathbb{Z}}
\def\NN{\mathbb{N}}

\begin{document}

\title{Singularity barriers and Borel plane analytic properties of $1^+$
  difference equations} \author{O.  Costin}\gdef\shorttitle{Borel plane
  analyticity of $1^+$ equations}

\maketitle
\date{}
\bigskip

\begin{abstract}

  The paper addresses generalized Borel summability of ``$1^+$'' difference
  equations in ``critical time''.  We show that the Borel transform $Y$ of a
  prototypical such equation is analytic and exponentially bounded for
  $\Re(p)<1$ but there is no analytic continuation from $0$ toward $+\infty$:
  the vertical line $\ell:=\{p:\Re(p)=1\}$ is a singularity barrier of $Y$.
  
  There is a unique natural continuation through the barrier, based on the
  Borel equation dual to the difference equation, and the functions thus
  obtained are analytic and decaying on the other side of the barrier.  In
  this sense, the Borel transforms are analytic and well behaved in
  $\CC\setminus\ell$.

  The continuation provided allows for generalized Borel summation of the
  formal solutions. It differs from standard ``pseudocontinuation''
  \cite{Shapiro}. This stresses the importance of the notion of cohesivity, a
  comprehensive extension of analyticity introduced and thoroughly analyzed by
  \'Ecalle.
  
  We also discuss how, in some cases, \'Ecalle acceleration can provide a
  procedure of natural continuation beyond a singularity barrier.

\end{abstract}

\section{Introduction} 

In the case of generic differential equations, generalized Borel summation of
a formal power series solution, in the sense of \'Ecalle \cite{Ecalle},
essentially consists in the following steps: (1) Borel transform with respect
to a {\em critical time}, related to the order of exponential growth of
possible solutions, (see also the note below), usual summation of the obtained
series, analytic continuation along the real line or in its neighborhood,
proper averaging of the analytic continuations (e.g.  medianization) toward
infinity, possible use of acceleration operators and Laplace transform
$\mathcal{L}$.

The choice of the critical time, or of a very slight perturbation --weak
acceleration-- of it is crucial for \'Ecalle summability.  A slower variable
(time) would hide the resurgent structure encapsulating the Stokes phenomena,
and, perhaps more importantly, introduces superexponential growth preventing
Laplace transformability at least in some directions. In a faster variable,
convergence of the Borel transformed series would not hold.

In some functional equations and so called type $1^+$ difference equations,
new difficulties occur. For them, \'Ecalle replaces analyticity with {\em
  cohesivity} \cite{Ecalle2}. This property was studied rigorously for some
classes of difference equations by Immink \cite{Immink}.  It is the purpose of
this note to show the importance of this notion: even in simple $1^+$
difference equations it is shown that critical time Borel transform has
barriers of singularities, preventing continuation in some half-plane.  This
occurs in the prototypical equation

\begin{equation}
   \label{f1} y(x+1)=\frac{1}{x}y(x)+\frac{1}{x} \end{equation} (example
 2. of \cite{Immink}). A simple proof of Borel space natural boundaries is not
 present in the literature, as far as the author is aware.  We also show that
 the barrier is traversable: on the real line the associated function is well
 defined and Laplace transformable to a solution of the difference equation.
 This function is real analytic except at one point and, in fact has analytic
 continuation in the whole of $\CC\setminus\ell$ with $\ell=\{p:\Re(p)=1\}$ a
 singularity barrier.  The present approach is adaptable to more general
 equations.
 
 We expect barriers of singularities to occur quite generally in $1^+$ cases,
 due to the fact that the pole position is periodic in the original variable,
 while critical time introduces a logarithmic shift in this periodicity. This
 leads to lacunary series in Borel plane, hence to singularity barriers.
 
 Nonetheless, further analysis shows that, in this simple case, and likely in
 quite some generality, softer Borel summation methods and study of Stokes
 phenomena are possible, relying on the convolution equation for continuation
 through singularity barriers.
 
 In spite of its simplicity, the properties in Borel plane of this equation,
 in the critical time, are very rich. 
 
 \z {\bf Note on critical time.}  The solution of the homogeneous equation
 associated to (\ref{f1}), $f(x)=1/\Gamma(x)$ has large $x$ behavior
 $(x/{2\pi})^{1/2}e^{-x\ln x+x}$.  The critical time $z$ is then the leading
 asymptotic term in the exponent, $z=x\ln x$ \cite{Immink}.  (The origin of
 the terminology $1^+$ is related to the exponential order slightly larger
 than one of $f$).  Various slight perturbations of this variable, weak
 accelerations, are used and indeed are quite useful.

\section{The singularity barrier} 
\begin{Theorem}\label{barr}
  Let $Y(p)$ be the Borel transform of $y$ in (\ref{f1}) in the critical time
  $z$.  Then $Y(p)$ is analytic on $\{p\ne 0:\arg(p) \in (\pi-2\pi,\pi+2\pi);
  \Re(p)<1\}$ and exponentially bounded as $|p|\to\infty$ in this region.  The
  line $\ell=\{p:\Re p=1\}$ is a singularity barrier of $Y$.
\end{Theorem}

{\em Proof of the theorem.}  Let $\tilde{y}$ be the formal power series
solution of (\ref{f1}). We study the analytic properties of the
Borel transform $\mathcal{B}\tilde{y}:=Y(p)$ of the on $\mathcal{S}_0$, the
Riemann surface of the log at zero, with respect to the critical time $z$.  In
critical time the functional equation of $\mathcal{B}\tilde{y}$
(\ref{eq:complicatedeq}) is unwieldy, and instead we look at the meromorphic
structure of solutions on which we perform a Mittag--Leffler decomposition.

It is straightforward to check that $\tilde{y}$ is the asymptotic series for
$\arg(x)\ne 0$ of the following actual solution of (\ref{f1})
\begin{equation}
  \label{e01}
  y_0(x)=\sum_{k=1}^{\infty}\prod_{j=1}^k\frac{1}{x-j}
\end{equation}
The fact that Res\,$(y_0;x=n)= e^{-1}/\Gamma(n)$ and the behavior at infinity
of $y_0$ show that the  Mittag-Leffler partial fraction decomposition of
(\ref{e01}) is
\begin{equation}
  \label{e01}
y_0=e^{-1}\sum_{k=1}^{\infty}\frac{1}{(x-k)\Gamma(k)}
\end{equation}

(1) {\em Analyticity in the left half plane.} The inverse function $z\mapsto x(z)$ of
$x\ln x$ is analytic on $\mathcal{S}_0\setminus (-e^{-1},0)$ as it can be seen
from the differential equation $\frac{dx}{dz}=(1+\ln x)^{-1}$.  Then $Y(p)$ is
the analytic continuation of the function defined for $p$ {\em negative} by
\begin{equation}
  \label{s2}
-\frac{1}{2\pi i}\int_{ i\RR-e^{-1}}e^{pz}y_0(x(z))dz=\frac{1}{2\pi i}\int_{C}e^{pz}y_0(x(z))dz,\ \  p\in\RR^-
\end{equation}
where $C$ is a contour  from $\infty+i0$ around
$ -e^{-1}$ and to $\infty-i0$.

(2) {\em Identities for finding continuation in  $\{z:\Re(z)<1\}$ and exponential bounds.}
 For analytic continuation clockwise we start from $\arg
p=\pi$ and rotate up the contour, collecting the residues:
\begin{eqnarray}
  \label{s4}
Y(p)=\frac{1}{2e\pi i}\sum_{k=1}^{\infty}\frac{1}{\Gamma(k)}\int_{C}\frac{e^{pz}dz}{x(z)-k}=F(p)+\frac{1}{2e\pi i}\int_{C_1}\sum_{k=1}^{\infty}\frac{1}{\Gamma(k)(x(z)-k)}e^{pz}dz\nonumber \\\ \text{where }\ F(p):=\sum_{k=1}^{\infty}\frac{1+\ln
  k}{e\Gamma(k)}e^{p k\ln k}
\end{eqnarray}
and where for small $\phi>0$, $C_1$ is the contour from $\infty e^{i\phi+i0}$
around $(-e^{-1},0)$ to $\infty e^{i\phi-i0}$. As $\arg p$ is decreased from
to zero (and further to $-\pi$), $\phi$ can be increased from $0^+$ to
$2\pi^-$ making $\int_{C_1}$ visibly analytic in $\{p\ne 0:\arg p
\in(-\pi,\pi)\}$ and exponentially bounded as $|p|\to \infty$. We decomposed
$Y$ into a sum of a lacunary Dirichlet series and a function analytic in the
right half plane. 

(2) {\em The natural boundary.} The Dirichlet series $F$ is manifestly
analytic for $\Re p <1$. As $p\uparrow 1$ we have $F(p)\to +\infty$ and thus
$F$ is not entire.  But then, by the Fabry-Wennberg-Szasz-Carlson-Landau
theorem \cite{Mandelbrojt} pp. 18, $\ell$ is a singularity barrier of $F$ and
thus of $Y$. For a detailed analysis, see also the note below. $\Box$
  
\bigskip {\bf Note: Description of the behavior of $F$ at $\ell$}. Since all
terms of the Dirichlet series are positive on the real line, it is easy to
check using discrete Laplace method\footnote{Determining, for fixed $p$, the
  maximal term of the series and doing stationary point expansion nearby.}
that $F$ increases like an iterated exponential along $\RR^+$ toward $\ell$,
$F(p)\propto \exp((1-p)\exp(1/(1-p)))$.  There are densely many points near
$\ell$ where the growth is similar; it suffices to take a sequence of
$k\in\NN$, $\Re(p)=k/(1+\ln(k))$ and $(1+\ln(k))\Im (p)$ very close to an
integer multiple of $2\pi$. (A Rouch\'e type argument shows there are also
infinitely many zeros with a mean separation of order the reciprocal of the
maximal order of growth, $\ln(d)\sim -(1-p)e^{1/(1-p)}$.)
Rather than attempting some form of continuation through points where $F$ is
bounded, which are easy to exhibit, we prefer to soften the barrier first, by
acceleration techniques.

\section{General Borel summability in the direction of the barrier. Properties
  beyond the barrier.}  

\z {\bf Strategy of the approach}.  It is convenient to perform a ``very weak
acceleration'' to smoothen the behavior of $Y(p)$ near $\ell$. The natural
choice of variable is $z=\ln\Gamma(k)$, but we prefer to slightly accelerate
further, to $z_m(x)$ defined in Remark~\ref{R11} below. We construct actual
solutions of (\ref{f1}) starting from an incomplete Borel sum. We identify
these actual solutions and show they are inverse Laplace transformable.
Furthermore, they solve the associated convolution equation in Borel space.
From these points of view, we have a unique continuation on $\RR^+$. We show
that the function thus obtained is real analytic on $\RR\setminus\{1\}$ and
continuable to the whole of $\CC\setminus \ell$.

The general solution of (\ref{f1}) is
\begin{equation}
  \label{eq:gensol}
  y(x)=y_0(x)+\frac{f(x)}{\Gamma(x)}
\end{equation}
where $f$ is any periodic function of period one, as it can be easily seen by
making a substitution of the form (\ref{eq:gensol}) in the equation. It can be
easily checked that the following solution of (\ref{f1})
\begin{equation}
  \label{eq:eqf2}
  y_1(x)=y_0+\frac{\pi}{ e}\frac{\cot\pi x}{\Gamma(x)}
\end{equation}
is an entire function, and has the asymptotic behavior $\tilde{y}$, the formal
series solution to (\ref{f1}) defined in the proof of the theorem. 
\begin{Remark}\label{R11}
  Let $m\in\NN$ and $z_m(x)=x\ln x -x-(m+\frac{1}{2})\ln x$. For given $C>0$,
  there is a one--parameter family of solutions of (\ref{f1}) which are
  analytic and polynomially bounded in a region of the form $S_C=\{x:\Re
  (z_m(x))\ge C\}$. They are of the form $y_c(x)=y_1(x)+c/\Gamma(x)$ for some
  constant $c$.
\end{Remark}

\z{\em Proof}. The solution (\ref{eq:eqf2}) already has the stated boundedness
and analyticity properties (and in fact, it decreases at least like $x^{-m}$
in $S_C$).  The general solution is of the form $y_1+f(x)/\Gamma(x)$ with $f$
periodic, as remarked at the beginning of the section. Analyticity implies $f$
is analytic and boundedness in the given region implies $f$ is bounded on the
line $\partial S_C$. By periodicity, $f$ is polynomially bounded in the whole
of $\CC$, which means $f$ is a polynomial, and by periodicity, a constant.
$\Box$

\begin{Theorem}[Generalized Borel summability]\label{BorelSum}
  (i) There exists a one parameter family of solutions of (\ref{f1}) which can
  be written as $\mathcal{L}_{z_m}H_c:=\int_0^{\infty}e^{-z_m p}H_c(p)dp$ where
    $H_c=\mathcal{B}_{z_m}\tilde{y}$ is analytic and exponentially bounded for
    $\Re(p)<1$ and $H_c\in C^{m-1}(\RR^+)$. 
  
    (ii) $H_c$ are real analytic on $\RR^+\setminus\{1\}$; they extend
    analytically to $\CC\setminus\ell$, and $\ell$ is a singularity barrier
    $H_c$ and the functions are $C^{m-1}$ on the two sides of the
    barrier\footnote{The values on the two sides cannot, obviously, be the
      same.}.  Furthermore, for $\Re(p)>1$, $H_c$ decrease toward infinity in
    $\CC$.
\end{Theorem}
\begin{Remark}
  It would not be correct at this time to conclude that, say,
  $\mathcal{L}^{-1}y_1$ provides Borel summation of $\tilde{y}$; we need to
  show that $y_1$ satisfies the necessary Gevrey--type estimates to identify
  the inverse Laplace transform with $\mathcal{B}\tilde{y}$ in the unit disk.
  We prefer to proceed in a more general way, not using explicit formulas, but
  constructing actual solutions starting with an incomplete Borel summation
  (and identifying them later with the explicit formulas).
\end{Remark}
{\em Proof of Theorem \ref{BorelSum}, (i)} We redo the analysis of the proof
of Theorem 1 in the variable $z=z_m$ and we get a decomposition of the form
(\ref{s4}), where now $F$ is replaced by
\begin{equation}
  \label{eq:s5}
  F_2=\sum_{k=1}^{\infty}\frac{\ln
  k+\frac{m}{k}}{e\Gamma(k)}e^{p \big[k\ln k-k-(m+\frac{1}{2})\ln k\big]}
\end{equation}
which is a Dirichlet series of the same type as $F$ and hence has $\ell$ as a
singularity barrier. However, $F_2$ is (manifestly) uniformly $C^{m-1}$ up to
$\ell$ and so is thus $Y(p)$.

For the solutions of (\ref{f1}) that decrease in a sector in the right half
--plane it is clear that the dominant balance is between $y(x+1)$ and $1/x$.
We then rewrite the equation to prepare it for a contraction mapping argument
in Borel space. By a slight abuse of notation we write $y(z)$ for $y(x(z))$
and we have
 $$(x(z)-1)y(x(z))=y(x(z)-1)+1$$
$$(x(z)-1)y(z)=y(z-g(z))+1$$
where $g(z)=\ln z-\ln\ln z+o(1)$ and then
$$(x(z)-1)y(z)=\sum_{k=0}^{\infty}y^{(k)}(z)g(z)^k/k!+1$$
Thus, dividing by
$x(z)-1$ and taking inverse Laplace transform,  with
$G_k(p)$ the inverse Laplace transform of ${g(z)^k}/{(x(z)-1)}/k!$,  we have
\begin{equation}
  \label{eq:complicatedeq}
  Y(p)=\sum_{k=0}^{\infty}[(-p)^kY]*G_k(p)+F(p)
\end{equation}
The term $G_k$ is (roughly) bounded by $|e^{-k(1-p)}|$, as can be seen by the
saddle point method applied to the inverse Laplace transform integral. It is
easy to check, using standard contraction mapping arguments (see e.g.
\cite{Duke}), that $Y$ is given by a convergent ramified expansion in the open
unit disk. This was to be expected from estimates of the divergence type of
the formal solutions of (\ref{f1}). However, given the estimates on the terms
of the convolution equation,  the equation, as written, cannot be
straightforwardly interpreted beyond $\Re(p)=1$, the threshold of convergence
of the ingredient series. It is however possible to write a meaningful global
equation by returning to the definition in terms of Laplace transform. We then
write
$$\mathcal{L}^{-1}y(z+g(z))=\frac{1}{2\pi
  i}\int_{c-i\infty}^{c+i\infty}dz e^{pz}\int_0^{\infty}dq e^{-q(z+g(z))}Y(q)=
\int_0^{\infty}H(p,q)Y(q)dq$$
where 
$$H(p,q)=\frac{1}{2\pi i}\int_{c-i\infty}^{c+i\infty}e^{(p-q)z-qg(z)}dz=
\frac{1}{2\pi i}\int_{c-i\infty}^{c+i\infty}e^{(p-q)z+q(\ln\ln
  z+...)}z^{-q}dz$$
which is well defined for $q>0$ and integrable at $q=0$;
the convolution equation becomes
\begin{equation}
  \label{eq:conveq}
  \int_0^{\infty}H(p,q)Y(q)dq=Y*\mathcal{L}^{-1}\Big[\frac{1}{x(z)-1}\Big]+
\mathcal{L}^{-1}\Big[\frac{1}{x(z)-1}\Big]
\end{equation}

Based on the solution on $[0,1)$ of (\ref{eq:complicatedeq}) we
  construct solutions to (\ref{f1}) and their inverse Laplace transforms
  provide continuation of $Y$ past $\Re(p)=1$ and implicitly solutions to
  (\ref{eq:conveq}).

 We define the incomplete Borel sum
 $$\hat{y}=\int_0^1 e^{-zp} Y_1(p)dp$$
 Formal manipulation shows that
 $\hat{y}$ satisfies (\ref{f1}) with errors of the form\footnote{Resulting
   from incomplete representation of $1/(x(z)-1)$.}  $o(e^{-z})$ or
 $o(x^{m}/\Gamma(x))$ in the variable $x$)  where the estimate of the errors
 is uniform in the right half--plane in $z$, or in a region $S_C$ w.r. to $x$.
 
 We look for a solution of (\ref{f1}) in the form
 $\hat{y}+\delta(x)/\Gamma(x)$.  Then $\delta(x)$ satisfies
 $\delta(x+1)=\delta(x)+R(x)$ (the $1^+$ degeneracy is not present anymore)
 where $R(x)=o(x^m)$ with differentiable asymptotics (by Watson's lemma).  A
 solution of this equation is
 $-\mathcal{P}^{m+3}\sum_{k=0}^{\infty}R^{(m+3)}(x+k)$, with $\mathcal{P}$ an
 antiderivative,  which is manifestly
 analytic and polynomially bounded  in regions of the form $S_C$, and
 $\hat{y}+\delta/\Gamma$ is manifestly a solution of (\ref{f1}), which, by
 construction, is also polynomially bounded  in $S_C$.
 
 By Remark \ref{R11},  $\hat{y}+\delta/\Gamma$ is one of the  solutions $y_c$. 
 But  $y_c$ is  inverse Laplace  transformable with  respect to  $z$,  and has
 sufficient  decay  to  ensure  the  existence of  $m-1$  derivatives  of  the
 transform. By  Remark \ref{R11}, any  solution that decreases in  the natural
 region $S_C$ in the right half plane  can be represented in this way and thus
 the conclusion follows.$\Box$
 \begin{Corollary}
   In $\{p:\Re(p)<1\}\cup [1,\infty)$, there is a one parameter family of
   Laplace transformable solutions to (\ref{eq:conveq}), the functions $H_c$
   in Theorem \ref{BorelSum} (i).  They have $\ell$ as a barrier of
   singularities.
 \end{Corollary}
\z {\em Proof of Theorem \ref{BorelSum} (ii)}. Since all Laplace transformable
solutions to (\ref{eq:conveq}) are those provided in Remark \ref{R11}, we
analyze the properties of the inverse Laplace transform of these functions for
$\Re(p)>1$. 

We note that, due to the fact that $y_c(z_m)$ increase at most as
$e^{z_m}/z_m^m$, we can deform for $\Re(p)>1$, the integral

\begin{equation}
  \label{eq:contdeform}
  \int_{c-i\infty}^{c+i\infty} e^{pz_m} y_c(z_m) dz_m
\end{equation}
to an integral
\begin{equation}
  \label{eq:contdeform}
  \int_C e^{pz_m} y_c(z_m) dz_m
\end{equation}
\z where $C$ starts at $-\infty-i\epsilon$, avoids the origin through the
right half plane and turns back to $-\infty+i\epsilon$. In view of the bound
mentioned above for $ y_c(z_m)$, this function is manifestly bounded and
analytic for $\Re(p)>1$, and in fact is continuous with $m-1$ derivatives up
to $\Re(p)=1$.  

\z {\bf Cohesive continuation and pseudocontinuation}.  It follows
from our analysis and from the fact that \'Ecalle's cohesive continuation also
provides solutions to the equation, that the results of the continuations are
the same (modulo the choice of one parameter, discussed in the Appendix). This
type of continuation is the natural one since it provides solutions to the
associated convolution equation. It is easy to see however that this
continuation is not a classical pseudocontinuation through the barrier, as it
follows from the following Proposition.

\begin{Proposition}
  The values of $H_c$ on the two sides of $\ell$ are not pseudocontinuations
  \cite{Shapiro} of each--other.
\end{Proposition}
\z {\em Proof.} Indeed, pseudocontinuation \cite{Shapiro}, pp. 49 requires that the
analytic elements coincide almost everywhere on the two sides of the barrier.
But $H_c$ is continuous on both sides, and then the values would coincide
everywhere, immediately implying analyticity through $\ell$, a contradiction.
\begin{Remark}
  The axis $\RR^+$, which is also a Stokes line, plays a special role. No
  other points on the singularity barrier can be used for Borel summation, as
  shown in the proposition below.
\end{Remark}
\begin{Proposition}\label{ExpGrowth}
  No Laplace transformable solution of (\ref{eq:conveq}) exists, in directions
  $e^{i\phi}\RR^+$, $\phi\in(0,\pi/2)$. (The same conclusion holds with
  $\phi\in(-\pi/2,0)$.)
 \end{Proposition}
 {\em Proof}. Indeed, the Laplace transform $y$ of such a solution would be
 analytic and decreasing in a half plane bisected by $e^{i\phi}$ and
 solve(\ref{f1}). Since $1/\Gamma(x)$ is entire and the general solution is of
 the form (\ref{eq:gensol}), by periodicity $f_1=f-\frac{\pi}{e}\cot\pi x$
 would be entire too.  Taking now a ray $te^{i(\phi+\pi/2-\epsilon)}$ we see,
 using again periodicity, that $f_1$ decreases factorially in the upper half
 plane.  Standard contour deformation shows that half of the Fourier
 coefficients are zero, $f_1 (x)=\sum_{k\in\NN}c_k e^{ikx}$ and that, because
 $f$ is entire, $c_k$ decrease faster than geometrically. But then
 $f_1(x)=:F(\exp(2\pi ix))$ with $t\mapsto F(t)$ entire. When $x\to i\infty$,
 $t\to 0$ and, unless $F=0$, we have $F(t)\sim ct^n$ for some $n\in\NN$, thus
 $f(x)\sim ce^{inx}$, incompatible with factorial decay.  This means $f=0$ but
 then (\ref{eq:gensol}) is not analytic on the real line\footnote{We should
   note that a procedure mimicking the proof of Theorem \ref{BorelSum} (i) in
   non-horizontal directions would fail because now the remainders $R(x)$
   would grow fast along the direction of evolution -- parallel to $\RR^+$.}.
 $\Box$

\section{Appendix: Weak acceleration, integral representation, median choice, natural
  crossing of the barrier}\label{Cont} A weak acceleration is provided by the
passage $x\ln x -x\mapsto x$. The $x$-- inverse Laplace transform of
(\ref{f1}) satisfies $e^{-p}Y-\int_0^p Y(s)ds -1=0$ with the solution
$Y=e^{-1}\exp(p+\exp(p))$.  $\mathcal{L}Y$ exists along any (combination of)
paths $R_n$ starting from the origin and ending on a ray of the form
$p=\RR^++(2n+1)i\pi, n\in\ZZ$.  The function
$f_+=\int_{R_1}e^{-xp}e^{p+e^p-1}dp$ is manifestly entire\footnote{It
  provides, in view of the superexponential properties of the integrand, Borel
  oversummation.}.  For $x=-t;t\to\infty$ the saddle point method gives
$$
f_+\sim \sqrt{2\pi}e^{t\ln t-t+\pi i t+\frac{1}{2}\ln t-1}$$
which
identifies $f_+$ with $ y_1+\pi i/e/\Gamma(x)$. With obvious notations, we see
that $y_1=\frac{1}{2}(f_++f_-)$, reminiscing of medianization. We have also
checked numerically that $y_1$ is approximated by least term truncation of its
asymptotic series with errors $o(1/\Gamma(x))$. (The integral representation
would allow for a rigorous check, but we have not done this and we state the
property as a conjecture; we also conjecture that the solution constructed in
Proposition \ref{BorelSum} is $y_1$; this could be checked by looking at the
asymptotic behavior on $\partial S_C$.) There is, obviously, only one solution
so well approximated. It should then be considered as the natural candidate
for the medianized transform in critical time and its inverse Laplace
transform, defined on the whole of $\RR^+$, and the natural continuation of
the Borel transform $\mathcal{B}\tilde{y}$ past the barrier.  For all these
reasons it is likely, but we have not checked it rigorously, that $y_1$
corresponds to the medianized cohesive continuation of \'Ecalle.
\begin{Remark}
  The procedure described of naturally crossing a barrier does not necessarily
  depend on the existence of an underlying functional equation.  It is
  sufficient to have accelerations as above that allow for Borel
  (over)summation along some paths, and choose as a natural actual function
  the one that has minimal errors in least term truncation or resort to a
  medianized choice. The process of continuation through the barrier can be
  written as the composition
  $\mathcal{L}^{-1}_{z_m}\mathcal{L}_{z_1}\mathcal{B}_{z_1}\hat{\mathcal{L}}_{z_m}$
  with $\hat{\mathcal{L}}$ formal Laplace transform, and is expected to
  commute with most operations of natural origin. It is applicable to many
  other series including the Dirichlet series $\sum _{k=0}^{\infty}
  e^{(p-1)n^2}$.
\end{Remark}

Finally, it seems a plausible conjecture that in the case of nonlinear
systems, infinitely many equally spaced ``isolated'' barriers should occur.

\z {\bf Acknowledgments}. The author is grateful to B. L. J.  Braaksma, and G.
Immink for pointing out to the problem and for very useful discussions and to
R. D.  Costin for a valuable technical suggestion. The work was partially
supported by NSF grant 0406193.

\end{document}